\documentclass[letterpaper, 10 pt, conference]{ieeeconf}
\usepackage{amsmath, amssymb,amsfonts}
\usepackage{subcaption}
\usepackage{graphicx}
\usepackage{color}
\usepackage{algorithm}
\usepackage{algpseudocode}
\usepackage{hyperref}

\usepackage{multirow}
\usepackage{changes}

\newtheorem{theorem}{Theorem}
\newtheorem{definition}{Definition}

\newtheorem{remark}{Remark}
\newtheorem{example}{Example}

\begin{document}

\title{\LARGE \bf
Stable MPC with maximal terminal sets and quadratic terminal costs
}

\author{Mikael Johansson and Hamed Taghavian\\
	 School of Electrical Engineering and Computer Science, KTH, Stockholm, Sweden}

\maketitle

\begin{abstract}
This paper develops a technique for computing a quadratic terminal cost for linear model predictive controllers that is valid for all states in the maximal control invariant set. This maximizes the set of recursively feasible states for the controller,  ensures asymptotic stability using standard proofs, and allows for easy tuning of the controller in linear operation.
\end{abstract}

\section{Introduction}

With its ability to optimize the closed-loop system under constraints on states and control signals, 
model predictive control (MPC) 
has emerged as the preferred control technique for many advanced applications. Although the underlying ideas can be traced back to the 1960's~\cite{Pro:63} and industrial applications appeared already in the 70's~\cite{Ric:77}, a more complete theoretical understanding of MPC emerged first around the millenium~\cite{May:00}. Advances in optimization algorithms~\cite{BEMPORAD20023,Fer:14,osqp} now allow model-predictive controllers to run on cheap embedded hardware, and a wealth of applications have demonstrated the practical value of MPC~\cite{RaL:18}. With a number of excellent textbooks on MPC~\cite{BBM:17,KoC:16,MRR+:00}, it is fair to say that the field of model predictive control for linear constrained systems is rather complete. However, there are still a few important questions that have not been fully resolved. This paper aims to address one of them.

Specifically, we are concerned with the problem of achieving the maximal operating region of a linear MPC with quadratic stage and terminal costs, while guaranteeing asymptotic stability of the closed loop.

The standard stability theorems for linear MPC require that the terminal set is matched with an appropriate terminal cost. Textbooks on MPC, such as~\cite{BBM:17,KoC:16,MRR+:00}, only present partial and suboptimal solutions to this problem. The easiest approach is to require that the terminal state reaches the origin, in which case we do not need a terminal penalty. A less restrictive solution is to use the control invariant set induced by \emph{some} linear state feedback controller as terminal set. A matching terminal cost can then be computed by solving a Lyapunov equation. However, both these approaches result in an unnecessarily small feasible region for the associated MPC controller. To increase their operating regions, one needs to use long prediction horizons. A more economical alternative is to use the maximal control invariant set as a terminal set. 
Although techniques for computing maximal control invariant sets are well-developed~\cite{Bla:99}, it is not obvious how to match these with an appropriate terminal cost.
This paper demonstrates how one can compute a quadratic terminal cost that is valid on the maximal control invariant set, hence allowing the largest possible operating region for the associated MPC controller while guaranteeing asymptotic stability of the closed-loop.
    
The challenge of determining  terminal sets and terminal costs that enlarge the operating region for MPC controllers has been approached by several authors. For example,~\cite{DSM+02} used the invariant set of a saturated linear controller (instead of linear) as terminal set,~\cite{BLA:14} computed a piecewise affine terminal cost for hybrid MPC controllers, and~\cite{GrP:13} and~\cite{ScD:15} suggested to use terminal costs derived from the Minkowski (control) Lyapunov function induced by the terminal set. Although these latter two approaches enable the MPC controller to operate from all states in the interior of the maximal control invariant set~\cite{ScD:15}, they have a few disadvantages. First, polyhedral Lyapunov functions can have many segments, and therefore be complex to represent and costly to use in the MPC computations. This issue is partially overcome in~\cite{GrP:13}, where the same inequalities are used to define the terminal set and the terminal cost. The second, and more important, disadvantage of the approaches in~\cite{ScD:15,GrP:13} is that the computed terminal cost only depends on the system dynamics and constraints, but is independent of the stage cost. This necessitates a different convergence proof and makes it hard to affect the properties of the controller when it regulates the system without activating any constraints. 

This paper is organized as follows. Section~\ref{sec:linear_mpc} reviews the linear MPC set-up and the standard stability theorems. Our approach to terminal cost computation is developed in  Section~\ref{sec:tcost_computation}, and evaluated numerically in Section~\ref{sec:evals}. Finally, Section~\ref{sec:conclusions} concludes the paper. 

 \paragraph*{Notation.} Our notation is largely standard. 
For a set $\mathcal{C}\subseteq \mathbb{R}^n$, we let $\partial \mathcal{C}$ denote its boundary, and define the scaled set $\lambda\mathcal{C}=\{ x=\lambda y \;\vert\; y\in \mathcal{C}\}$. We say that $\mathcal{C}$ is a C-set if it is convex, compact, and includes the origin in its interior. 

A simplex ${\mathcal S}\subset \mathbb{R}^n$ is the convex hull of $n+1$ affinely independent vectors $v_i\in {\mathbb R}^n$ which we refer to as vertices. A triangulation of a set ${\mathcal C} \in {\mathbb R}^n$ is a subdivision of the set into a finite number of $n$-dimensional simplices such that any two simplices intersect in a common face (a simplex of any lower dimension) or not at all. In a boundary triangulation, each simplex in the triangulation has one vertex at the origin. 

\section{Linear Model Predictive Control} \label{sec:linear_mpc}
Consider the discrete-time linear system
\begin{align}
    x_{t+1} &= Ax_t + Bu_t,\quad t\geq 0
    \label{eqn:linsys}
\end{align}
with linear constraints on the states and controls
\begin{align}
    x_t\in \mathcal{X}, \qquad u_t\in \mathcal{U}. \label{eqn:csets}
\end{align}
Both $\mathcal{X}$ and $\mathcal{U}$ are assumed to be polyhedral C-sets.

It is useful to view linear MPC as an approximate solution to the infinite-horizon control problem
\begin{align}
    \begin{array}[c]{ll}
    \mbox{minimize} & \sum_{t=0}^{\infty} x_t^{\top}Qx_t + u_t^{\top}R u_t\\
    \mbox{subject to} & x_{t+1}=Ax_t + Bu_t\\
    & x_t\in \mathcal{X},\; u_t\in \mathcal{U}
    \end{array}\label{eqn:cstr_opt}
\end{align}
At every sampling instant $t$, the MPC control law measures the state $x_t$, solves the planning problem
\begin{align}
    \begin{array}[c]{ll}
    \underset{\{\hat{x}_k\}, \{\hat{u}_k\}}{\mbox{minimize}} & \sum_{k=0}^{T-1}\hat{x}_k^{\top}Q\hat{x}_k + \hat{u}_k^{\top}R\hat{u}_k + \hat{v}(\hat{x}_T)\\
    \mbox{subject to} & \hat{x}_{k+1} = A\hat{x}_k + B\hat{u}_k,\;\; k=0,\dots, T-1 \\
    & \hat{x}_k\in \mathcal{X},\; \hat{u}_k\in \mathcal{U},\;\;\;\;\;\;\; k=0, \dots, T-1\\
    &\hat{x}_T\in \mathcal{X}_T\\
    &\hat{x}_0 = x_t
    \end{array} \label{eqn:planning_problem}
\end{align}
for the predicted optimal controls $\{\hat{u}_k^{\star}\}$ and predicted state trajectory $\{\hat{x}_k^{\star}\}$, and
 applies the control 
\begin{align}
    u_t &= \hat{u}_0^{\star}. \label{eqn:rhc_control}
\end{align}
In the planning problem,  $\hat{v}(\hat{x}_T)=\hat{x}_T^{\top}Q_T\hat{x}_T$ serves as an approximation of the infinite-horizon cost-to-go of (\ref{eqn:cstr_opt}) from state $\hat{x}_T$ at the end of the planning horizon, while the terminal set $\mathcal{X}_T$ ensures that the cost-to-go approximation is valid. The standard stability proof for linear MPC imposes the following requirements \cite{BBM:17}.

\begin{theorem} \label{thm:mpc_stability}
Consider the system (\ref{eqn:linsys}) with constraints (\ref{eqn:csets}) under the the RHC control law (\ref{eqn:planning_problem})-(\ref{eqn:rhc_control}). Assume that $(A,B)$ is a reachable pair and let $\mathcal{X}_0$ be the set of states $x_t$ for which the planning problem (\ref{eqn:planning_problem}) admits a feasible solution. If
\begin{itemize}\setlength{\itemindent}{1em}
    \item[(a)] $Q\succeq 0$ with $(Q^{1/2},A)$ detectable, $R\succ 0$, $Q_T\succ 0$
    \item[(b)] The sets $\mathcal{X}$, $\mathcal{U}$ and $\mathcal{X}_T\subseteq X$ are C-sets
    \item[(c)] For every $x\in \mathcal{X}_T$, there exists a $u \in \mathcal{U}$ such that   
    \begin{align*}
        Ax+Bu \in \mathcal{X}_T,\, \mbox{ and }\\
        \hat{v}(Ax+Bu) - \hat{v}(x) + x^{\top}Qx + u^{\top}Ru \leq 0
    \end{align*}
\end{itemize}
Then, every trajectory $\{x_t\}$ of the closed-loop system remains in $\mathcal{X}_0$ and $\lim_{t\rightarrow\infty} x_t=0$.
\end{theorem}

While the first two conditions of the theorem are straightforward to verify, the last one is more involved. In essence, it requires that ${\mathcal X}_T$ is control invariant and that $\hat{v}(\cdot)$ is an upper bound on the true cost-to-go of (\ref{eqn:cstr_opt}) for all $x\in\mathcal{X}_T$. MPC textbooks, such as \cite{MRR+:00,KoC:16,BBM:17}, typically suggest two approaches to this problem. The first one is to set $\mathcal{X}_T=\{0\}$ (forcing the state at the end of the planning horizon to zero), for which we can set $\hat{v}_T(x)\equiv 0$. The second one is to choose a terminal set that is invariant under \emph{some} linear state feedback $u_t=-Lx_t$. Specifically, one uses the largest invariant set of $x_{t+1}=(A-BL)x_t$ contained in the set $\{ x\;\vert\; x\in \mathcal{X} \mbox{ and } -Lx\in \mathcal{U}\}$. A valid quadratic upper bound of the cost-to-go is then  $\hat{v}(x)=x^{\top}Q_Tx$ where $Q_T$ solves the Lyapunov equation
\begin{align*}
    Q_T &= Q + L^{\top}RL + (A-BL)^{\top}Q_T(A-BL) 
\end{align*}
A particularly convenient choice is to use the infinite-horizon optimal LQR controller $u_t=-L_{\infty}x_t$ for the cost defined by $Q$ and $R$. In this way, the Lyapunov equation above is satisfied for the solution $P_{\infty}$ to the corresponding discrete-time algebraic Riccati equation
\begin{align*}
    P_{\infty} &= Q+A^{\top}P_{\infty}A - L_{\infty}^{\top}(B^{\top}P_{\infty}B + R)L_{\infty}
\end{align*}
where $
    L_{\infty} = (B^{\top}P_{\infty}B + R)^{-1}B^{\top}P_{\infty}A$.
However, the terminal set computed in this way is not necessarily large.

Since $x_T$ must belong to $\mathcal{X}_T$,  a smaller terminal set leads to a smaller set of (recursively) feasible states, and hence to a smaller operating region of the MPC controller. The operating region increases with the planning horizon $T$, but with a small terminal set one typically need a long horizon to be able to operate from all states in the maximal control invariant set~\cite{ScR:94}. If one is able to use the maximal control invariant set as terminal set, on the other hand, the MPC controller will have the largest possible operating region already with a horizon of one. To illustrate the relationship between terminal set, prediction horizon, and operating regime of linear MPC, we consider the following example from~\cite{KoC:16}.

\begin{example} \label{ex:horizon}
Consider the second-order system
\begin{align*}
    x_{t+1} &= \begin{pmatrix}
    1.1 & 2\\ 0 & 0.95
    \end{pmatrix}x_t + 
    \begin{pmatrix}
        0\\0.0787
    \end{pmatrix}u_t \\
    y_t &= \begin{pmatrix}
        -1 & 1
    \end{pmatrix}x_t
\end{align*}
under the MPC control (\ref{eqn:planning_problem})-(\ref{eqn:rhc_control}) with
\begin{align*}
    Q &= C^{\top}C, \; R=1\\
    \mathcal{X} &= \{ x \;\vert\; \Vert x \Vert_{\infty}\leq 8\}\\
    \mathcal{U} &= \{ u \;\vert\; \vert u\vert\leq 1\}
\end{align*}
Figure~\ref{fig:invsets} shows $\mathcal{X}_0$ for different horizon lengths when the terminal set is taken to be the maximal invariant set of the infinite-horizon LQR-optimal control law. Note how the operating region of the MPC controller increases with increasing horizon length, and that even with a horizon of 24 samples, the MPC controller can only operate in a subset of the maximal control invariant set.
\begin{figure}[htb]
\centerline{\includegraphics[width=\hsize]{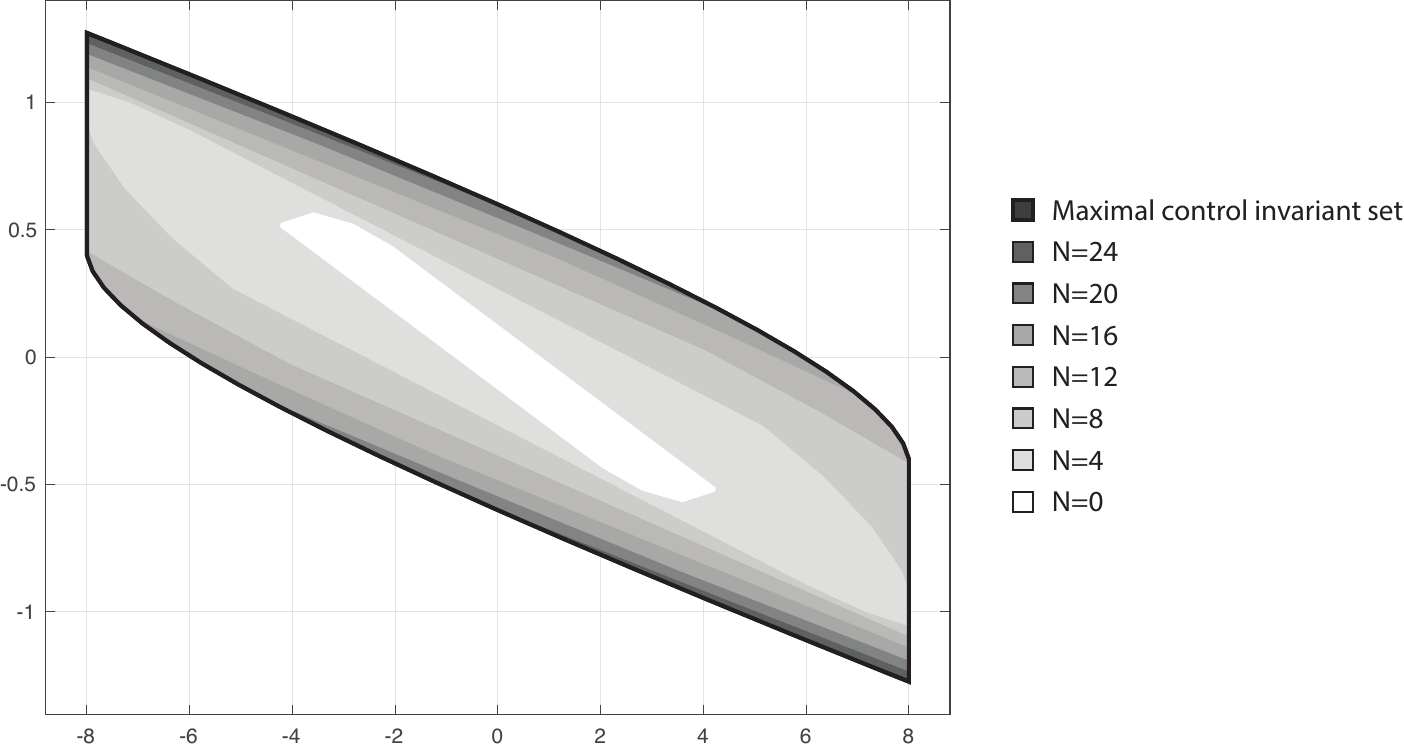}}
\caption{Operating region $\mathcal{X}_0$ for MPC controller depending on the horizon length $T$. A small terminal set may force us to use an unnecessarily long horizon.}
\label{fig:invsets}
\end{figure}
\end{example}

In the next section, we will develop a technique for computing a quadratic terminal cost so that one can use the maximal control invariant set as terminal set and still guarantee stability using Theorem~\ref{thm:mpc_stability}.

\section{Computing a quadratic terminal cost valid on the maximal invariant set}
\label{sec:tcost_computation}
Our procedure contains three key steps: we first determine the maximal invariant set, then recover an explicit feedback policy that renders the set invariant, and finally compute an upper bound on the infinite-horizon cost for this control law. 

\subsection{Maximal control invariant and $\lambda$-contractive sets}  \label{sec:cinv_sets}
\begin{definition} A set
$\mathcal{C}\subseteq \mathbb{R}^n$ is \emph{control invariant} for the system (\ref{eqn:linsys}) under the constraints (\ref{eqn:csets}) if $\mathcal{C}\subseteq\mathcal{X}$ and
\begin{align*}
    x &\in {\mathcal C}\Rightarrow \exists u\in \mathcal{U} \mbox{ such that } Ax+Bu\in \mathcal{C}.
\end{align*}
Furthermore, we say that $\mathcal{C}_{\infty}$ is maximal control invariant if it is control invariant and contains all control invariant sets for the system  (\ref{eqn:linsys}) under the constraints (\ref{eqn:csets}).
\end{definition}

Control invariance by itself does not imply that states in ${\mathcal C}_{\infty}$ can be driven to zero. For example, if $A=I$ and $B=0$, then ${\mathcal C}_{\infty}={\mathcal X}$ but $x_k=x_0$ for all $k\geq 0$. Such situations can be detected and avoided by requiring that the terminal set is contractive in the following sense.
\begin{definition} \label{def:lambda_contractive}
Let $\lambda\in (0,1]$. A set ${\mathcal C}\subseteq {\mathcal X}$ is called $\lambda$-contractive 
(\ref{eqn:linsys}) under the constraints (\ref{eqn:csets}), if for every $x\in {\mathcal C}$ there exists $u\in {\mathcal U}$ such that $Ax+Bu \in \lambda {\mathcal C}$. For a given $\lambda$, the maximal $\lambda$-contractive set for (\ref{eqn:linsys}) under the constraints (\ref{eqn:csets}), denoted ${\mathcal C}_{\infty}^{\lambda}$, is the union of all $\lambda$-contractive sets.
\end{definition}

Note that the maximal $1$-contractive set is identical to the maximal control invariant set. There are several techniques for computing the maximal control invariant set (\emph{e.g.},~\cite[Algorithm~10.2]{BBM:17}), but $\lambda$-contractive sets 
are typically computed using the following approach: define the predecessor set
\begin{align*}
    \mathcal{P}(\Omega) &= \left\{ x \in \mathbb{R}^n \;\vert\; Ax+Bu\in \Omega \mbox{ for some } u\in \mathcal{U} \right\} 
\end{align*}
and then execute the following algorithm~\cite[Algorithm~10.2]{BBM:17} 
\begin{enumerate}
    \item Set $\Omega_0=\mathcal{X}$,\quad $k=0$. 
    \item Let $k=k+1$ and $\Omega_{k}=\mathcal{P}(\lambda \Omega_{k-1}) \cap{\mathcal X}$.
    \item If $\Omega_{k}=\Omega_{k-1}$, return $C^{\lambda}_{\infty}=\Omega_k$ else goto 2.
\end{enumerate}

For a linear system with C-set constraints on the states and controls, the maximal control invariant set is also a C-set. However, even if $\mathcal{X}$ and $\mathcal{U}$ are polyhedral, the maximal invariant set may not be polyhedral and the algorithm above will not terminate. Conditions for when the algorithm is guaranteed to terminate in a finite number of steps, along with bounds for the number of iterations required, can be found in \cite{GiT:91}. When we consider $\lambda$-contractive sets with $\lambda\in [0,1)$, the situation becomes slightly more complicated since the constraints may limit our ability to contract; see~\cite{Bla:94} and ~\cite{SDC:17} for a careful convergence analysis of the algorithm. 

\subsection{Recovering the implicit feedback policy} \label{sec:cinv_feedback_policy}

The algorithm described above allows us to compute the maximal $\lambda$-contractive set $\mathcal{C}_{\infty}^{\lambda}$ of a linear system with C-set constraints on states and controls. From Definition~\ref{def:lambda_contractive}, we therefore know that there exists an admissible control function $u(x): \mathcal{C}_{\infty}^{\lambda}\mapsto \mathcal{U}$ that steers states in $\mathcal{C}_{\infty}^{\lambda}$ into $\lambda \mathcal{C}_{\infty}^{\lambda}$. However, this control function is only implicit in the calculations. We will now demonstrate how we can recover $u(x)$ as a continuous and piecewise linear control law. 

The next result, which is a slight generalization of the conditions proposed by Gutman and Cwikel~\cite{GuC:86} (see also~\cite[\S4]{Bla:99}) demonstrates how we can determine control signals to apply at each vertex of a C-set ${\mathcal C}$ to steer the state into $\lambda {\mathcal C}$.
\begin{theorem} \label{thm:vertex_controls}
The C-set $\mathcal{C}\subseteq\mathcal{X}$ with vertices $v_i$, $i=1, 2, \dots, s$ is $\lambda$-contractive for the discrete-time linear system (\ref{eqn:linsys}) under the constraints (\ref{eqn:csets}) if and only if there exists $\lambda_i\in [0,\lambda]$, $p_{ij}\in [0, 1]$ and $u_i\in \mathbb{R}^{m}$ that satisfy
\begin{align}\left\{ 
\begin{array}[c]{rcl}
    Av_i + Bu_i &=& \sum_{j=1}^s p_{ij} v_j \\
    \sum_{j=1}^s p_{ij} &\leq& \lambda_i  \\
    u_i &\in& \mathcal{U} 
\end{array}
    \right.
    \label{eqn:lp_recovery}
\end{align}
for every $i=1, 2, \dots, s$.
\end{theorem}

Together, the equality constraint 
and the conditions that  $p_{ij}\geq 0$ and $\sum_{j=1}^s p_{ij}\leq \lambda_i\leq \lambda$ imply that for each vertex $v_i$ of $\mathcal{C}$, there is an admissible control action $u_i$ so that the next state belongs to $\lambda_i {\mathcal C}$ and hence to $\lambda \mathcal{C}$. 
Since $\mathcal{U}$ is a C-set, $u\in\mathcal{U}$ can be expressed as a system of linear inequalities or equalities. This means that the conditions in (\ref{eqn:lp_recovery}) constitute a linear programming feasibility problem. 

\begin{remark}
Although we could fix $\lambda_i=\lambda$ for all $i$, it can be useful to minimize $\sum_i\lambda_i$ to encourage the control to drive the state closer to the origin or to minimize the difference between the $u_i$ and a linear control law, \emph{e.g.} $-L_{\infty}v_i$.
\end{remark}

As suggested by Gutman and Cwikel~\cite{GuC:86}, it is possible to transform the vertex controls computed in Theorem~\ref{thm:vertex_controls} into a continuous feedback policy that is valid for all $x$. To this end, consider a subdivision of ${\mathcal C}$ into simplices $\{ \mathcal{S}_k \}_{k=1}^{N}$, such that each simplex $\mathcal{S}_k$ has one vertex at the origin and the remaining $n$ ones at extreme points of $\mathcal{C}$. Such a boundary triangulation of a C-set is readily determined using standard convex hull algorithms \cite[Section~3]{BEF:00}. Let $v_i$ for $i=1, 2, \dots s$ be the extreme points of $\mathcal{C}$, $v_0=0$ and $\mathcal{I}_k = \{ i \;\vert\; v_i\in \mathcal{S}_k\}$ be the index set for the vertices of $\mathcal{S}_k$. 

Since ${\mathcal S}_k$ is a simplex, any $x\in  \mathcal{S}_k$ can be written as
\begin{align}
    x &= \sum_{i\in{\mathcal I}_k} p_i(x)  v_i 
    \label{eqn:x_vrep}
\end{align}
where $p_i(x)\geq 0$ and $\sum_{i\in {\mathcal I}_k} p_i(x) = 1$. Moreover, the vertices $v_i$ are affinely independent, so if we define 
$V_k\in \mathbb{R}^{n\times n}$ as the matrix whose columns are the non-zero vertices of ${\mathcal S}_k$ ordered in increasing vertex index $i$, the matrix
\begin{align*}
    \begin{pmatrix}
        0 & V_k \\
        1 & \mathbf{1}^{\top}
    \end{pmatrix}
\end{align*}
has full rank. This implies that $V_k$ also has full rank.

Next, define $u_0=0$ so that $Av_0+Bu_0=0$. A feasible solution to the conditions in Theorem~\ref{thm:vertex_controls} then implies that 
\begin{align*}
Av_i+Bu_i\in {\mathcal C}\qquad \mbox{ for all } i=0, 1, \dots, s.
\end{align*}
Since $\mathcal{C}$ is convex, any convex combination of $A v_i + Bu_i$ also belongs to ${\mathcal C}$. In particular, with $p_i(x)$ defined above,  
\begin{align*}
    \sum_{i\in \mathcal{I}_k} p_i(x) (Av_i + Bu_i) =Ax + B\sum_{i\in \mathcal{I}_k} p_i(x) u_i \in {\mathcal C}
\end{align*}
Similarly, since each $u_i \in {\mathcal U}$, and the set ${\mathcal U}$ is convex, 
\begin{align*}
    u(x) &= \sum_{i\in \mathcal{I}_k} p_i(x) u_i\, \in {\mathcal U}.
\end{align*}
Thus, $u(x)$ is admissible and renders ${\mathcal C}$ $\lambda$-contractive under $x_{k+1}=Ax_k+Bu(x_k)$. To derive an explicit expression for the control policy, let $U_k\in \mathbb{R}^{m\times n}$ be the matrix whose columns are $u_i$ for $i\in {\mathcal I}_k\backslash 0$ ordered in increasing vertex index $i$ and $p(x) \in \mathbb{R}^n$ as the vector of $p_i(x)$ for $i\in {\mathcal I}_k\backslash 0$ ordered in the same way. 
Then, for $x\in {\mathcal S}_k$, 
\begin{align*}
    x &= \sum_{i\in {\mathcal I}_k} v_i p_i(x) = V_k p(x) \Rightarrow p(x)=V_k^{-1}x
    \intertext{and}
    u(x) &= \sum_{i\in {\mathcal I}_k} u_i p_i(x) = U_k p(x)= U_k V_k^{-1}x 
\end{align*}
In other words, the feedback policy
\begin{align}
    u(x) &= -L_k x = U_k V_k^{-1}x \qquad x\in {\mathcal S}_k \label{eqn:pwl_sfb}
\end{align}
renders $\mathcal{C}$ $\lambda$-contractive. It is continuous and piecewise linear, and easy to extract from a solution to (\ref{eqn:lp_recovery}). 

\begin{example} \label{ex:trisets}
Figure~\ref{fig:trisets} shows the triangulation of the maximal control invariant set for the system considered in Example~\ref{ex:horizon}. The triangulation of $\mathcal{C}_{\infty}$ is shown in gray while the associated piecewise linear control law is shown in blue.
\begin{figure}[htb]
\centerline{\includegraphics[width=0.95\hsize]{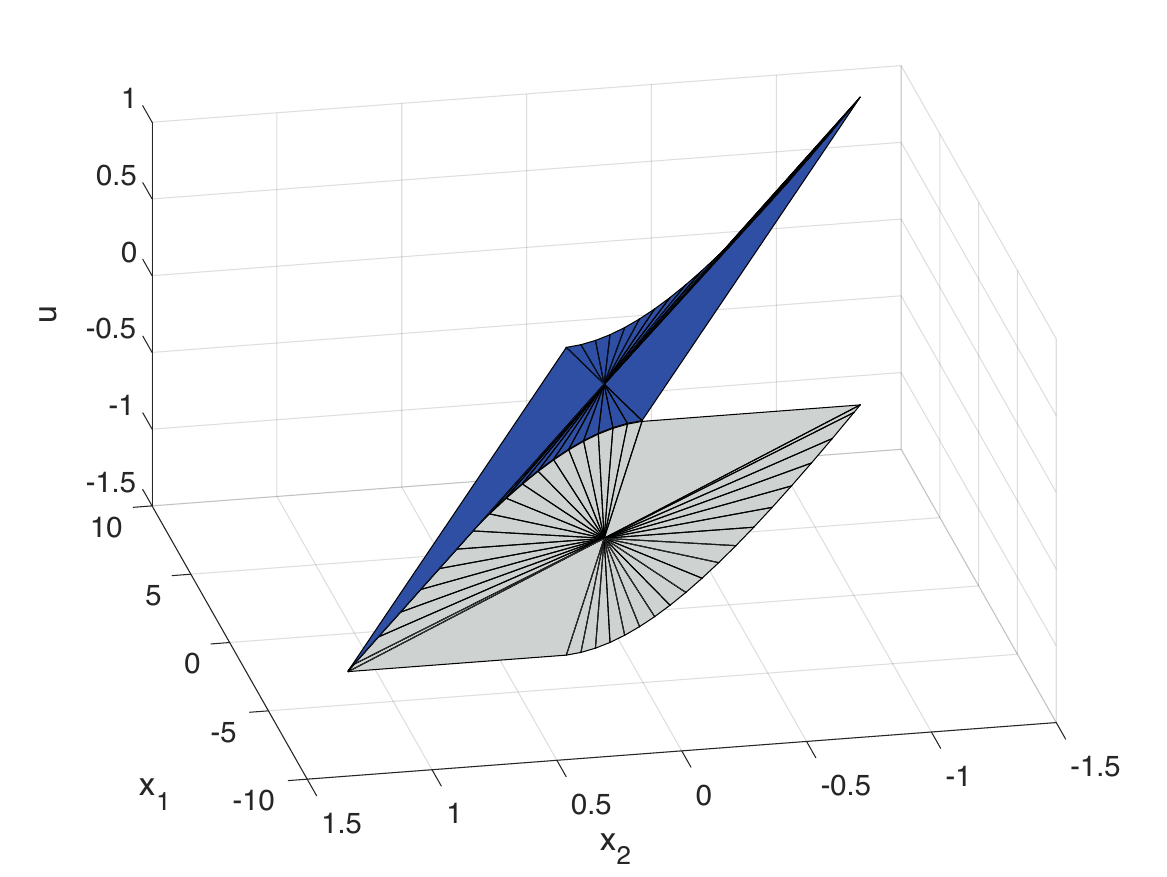}}
\caption{Triangulation of $\mathcal{C}_{\infty}$ (gray) and the associated piecewise linear control law (blue) for the system in Example~\ref{ex:horizon}.}
\label{fig:trisets}
\end{figure}
\end{example}

\subsection{A quadratic upper bound on the cost-to-go}
\label{sec:quad_bound}

Now that we have developed a procedure to recover an explicit feedback policy that renders the desired terminal set control invariant, we can proceed to compute an upper bound on the associated cost-to-go. In particular, 
consider $\hat{v}(x)=x^{\top}Px$. Then condition (d) of Theorem~\ref{thm:mpc_stability} reads
\begin{align*}
    x^{\top}A_k^{\top} P A_k x - x^{\top}Px + x^{\top}Q_k x &\leq 0 \quad \forall x\in {\mathcal S}_k
\end{align*}
where
$A_k=(A-BL_k)$ and $Q_k=Q+L_k^{\top}RL_k$. There are $N$ such inequalities, one for each simplex in the boundary triangulation of the terminal set. As described in~\cite{JoR:98,Joh:02} these conditions can be verified using semi-definite programming, where the condition that $x\in {\mathcal S}_k$ is encoded using the S-procedure. This leads to the next result.

\begin{theorem} \label{thm:quad_upper_bound}
Let ${\mathcal C}$ be control invariant for (\ref{eqn:linsys}) under constraints (\ref{eqn:csets}), and let $\{ \mathcal{S}_k\}_{k=1}^N$  be a boundary triangulation of $\mathcal{C}$. Consider the control policy $u_k=u(x_k)$ defined in (\ref{eqn:pwl_sfb}) and let $A_k=A-BL_k$ and $Q_k=Q+L_k^{\top}RL_k$. If the following semi-definite program
\begin{align*}
    \begin{array}[c]{lll}
    \mbox{minimize} & \mathrm{Tr}(P) &\\
    \mbox{subject to} & A_k^{\top}P A_k - P + Q_k + V_k^{-\top} W_k V_k^{-1} \preceq 0 & \forall k \\
    & W_k=W_k^{\top}\geq 0 & \forall k\\
    & \; P\succ 0 &
    \end{array}
\end{align*}
in variables $P$ and $W_k$ is feasible, then condition (d) in Theorem~\ref{thm:mpc_stability} is satisfied with $\mathcal{X}_T={\mathcal C}$ and $\hat{v}(x)=x^{\top}Px$.
\end{theorem}

\begin{remark}
The first matrix inequality uses the S-procedure to reduce conservatism. Since $W_k$ is elementwise non-negative, and $V_k^{-1}x\geq 0$ when $x\in \mathcal{S}_k$, the term $x^{\top}V_k^{-\top}W_kV_k^{-1}x\geq 0$ for $x\in \mathcal{S}_k$. However,  $V_k^{-\top}W_kV_k^{-1}$ is not necessarily a positive semi-definite matrix and $x^{\top}V_k^{-\top}W_kV_k^{-1}x$ may be negative when $x\not\in \mathcal{S}_k$. This makes the linear matrix inequality easier to satisfy than $A_k^{\top}PA_k-P+Q_k\preceq 0$. The S-procedure is, in general, only a sufficient condition for checking positivity of a quadratic form on a polyhedron, but since $\mathcal{S}_k$ is a simplex, this use of the S-procedure is loss-less up until dimension $n=4$~\cite{Joh:02}.     
\end{remark}

\begin{remark}
As discussed in \S~\ref{sec:cinv_feedback_policy}, the matrices $V_k$ are invertible per definition. Nevertheless, one can multiply the linear matrix inequalities
\begin{align*}
A_k^{\top}P A_k - P + Q_k + V_k^{-\top} W_k V_k^{-1} \preceq 0 
\end{align*}
by $V_k^{\top}$ from the left and $V_k$ from the right without affecting the solution set. By accounting for the structure of $A_k$, $Q_k$ and $L_k$, the resulting linear matrix inequality reads
\begin{align*}
(AV_k+BU_k)^{\top}P(AV_k+BU_k) + V_k^{\top}(-P+Q)V_k + \\U_k^{\top}RU_k + W_k \preceq 0
\end{align*}
This circumvents the need for recovering the feedback gains $L_i$ and avoids any inversions of the $V_k$ matrices.
\end{remark}

\begin{remark}
As a final remark, we have used the objective $\mbox{Tr}(P)$, but there are many other options. One such example would be to minimize $\sum_i v_i^TPv_i=\sum_i \mbox{Tr}\, P (v_i v_i^{\top})$, or to minimize $\Vert P-P_{\infty}\Vert_F^2$. Both are easily expressed as SDPs.
\end{remark}

\subsection{The complete algorithm}  \label{sec:algorithm}

We now have all the required components for computing a quadratic upper bound to the cost-to-go that is valid in the maximal control invariant set (and, more generally, in the maximal $\lambda$-contractive set). Specifically, we propose to use the following algorithm with $\lambda=1$.
\begin{enumerate}
    \item Compute  $\mathcal{C}^{\lambda}_{\infty}$ and determine its boundary triangulation (including the vertices) using a convex hull algorithm.
    \item Recover admissible vertex controls  $\{u_i\}$ that render ${\mathcal C}^{\lambda}_{\infty}$ $\lambda$-contractive by solving the linear program (\ref{eqn:lp_recovery}). 
    \item Convert $\{u_i\}$ into feedback gains $\{L_i\}$ such that
    \begin{align*}
        u_t &=-L_k x_t\quad \mbox{ for } x_t\in \mathcal{S}_k, \; k=1, \dots, N
    \end{align*}
    using the procedure in Section~\ref{sec:cinv_feedback_policy}.
    \item Solve the semidefinite program in Theorem~\ref{thm:quad_upper_bound} for $P$. 
\end{enumerate}
If the algorithm returns a feasible solution $P$, then the receding-horizon control law (\ref{eqn:planning_problem})-(\ref{eqn:rhc_control}) with $\mathcal{X}_T={\mathcal C}_{\infty}^{\lambda}$ and $Q_T=P$ renders the closed-loop control of the linear system (\ref{eqn:linsys}) asymptotically stable according to Theorem~\ref{thm:mpc_stability}.

\begin{remark}
Although the proposed procedure is numerical, the new steps that we have introduced rely on convex optimization, are fast and reliable to execute, and do not introduce any conservatism for systems of order $n\leq 4$. Still, since control invariance does not imply asymptotic stabilizability (cf. the discussion in \S~\ref{sec:cinv_sets}), we cannot guarantee that the procedure will always find a quadratic upper bound on the infinite-horizon cost-to-go that is valid in the maximal control invariant set. However, with $\lambda<1$, the implicit feedback policy is guaranteed to make the closed loop asymptotically stable and should, in principle, admit a quadratic upper bound on the infinite-horizon cost within the associated maximal $\lambda$-contractive set.
\end{remark}

\subsection{Horizon length and linear performance}
With the computed terminal cost, the MPC controller will result in an asymptotically stable closed-loop for all horizon lengths. Close to the origin, it will realize the linear state-feedback that is optimal for a $T$-horizon linear-quadratic control problem with stage cost defined by $Q$ and $R$, and terminal cost $x^{\top}Px$. The fact that the optimal unconstrained control is linear and easy to compute is a distinct advantage over approaches that use more complex terminal costs. It allows us to analyze the frequency domain properties of the controller near the origin, and understand how the horizon length affects the control performance. In particular, if the computed terminal cost matrix $P$ is very different from $P_{\infty}$ and we use a small horizon length, the linear operation of the MPC controller can be far from the infinite-horizon LQR controller. The difference diminishes with increasing horizon length, and is easy to quantify using the Riccati recursion for the associated finite-horizon LQR problem.
\section{Numerical examples}
\label{sec:evals}

In this section, we will illustrate various aspects of our framework by examples.  

\subsection{Horizon length and performance}

\begin{figure}[htb]
\centerline{\includegraphics[width=0.9\hsize]{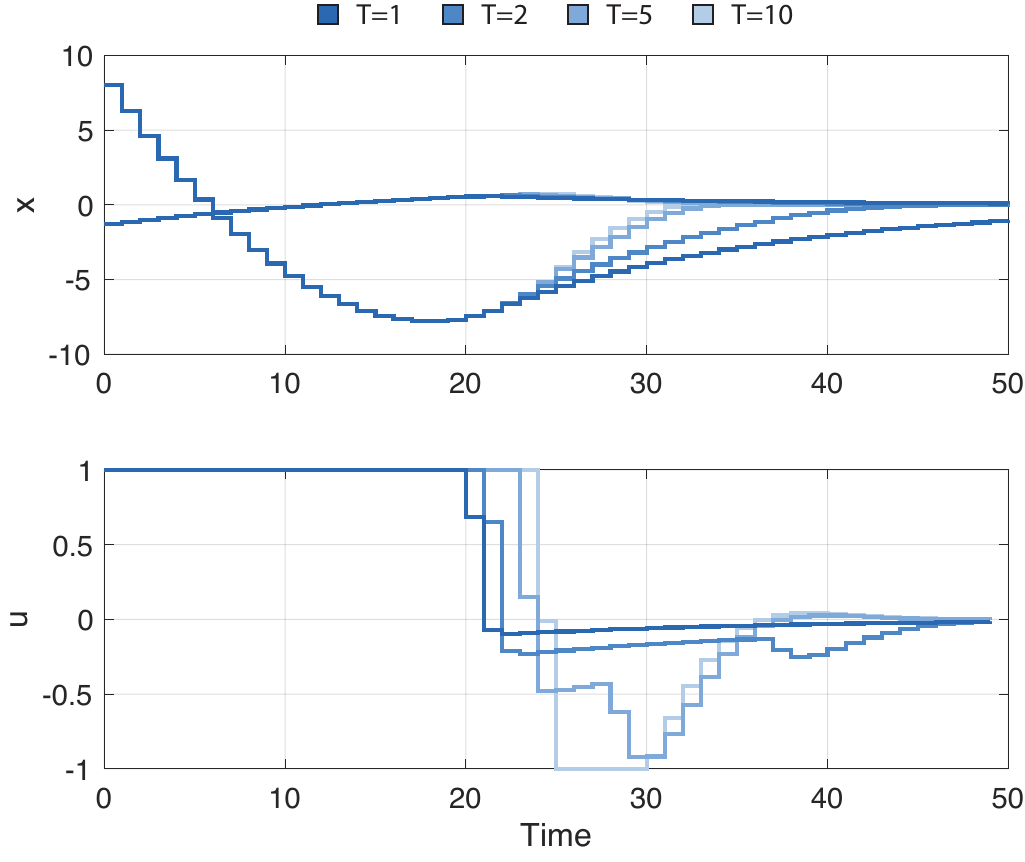}}
\caption{The MPC controller with the proposed terminal cost ensures a feasible and stable operation for all horizon lengths.}
\label{fig:sim2d}
\end{figure}

Let us first return to the system in Example~\ref{ex:horizon}. The maximal control invariant set ${\mathcal C}_{\infty}$ has $s=39$ vertices our procedure finds the terminal cost matrix
\begin{align*}
    Q_T &= \begin{pmatrix}
        38.6 & 343.1\\
        343.1 & 4178.5
    \end{pmatrix}
\end{align*}
If we compare this with the Riccati solution 
\begin{align*}
    P_{\infty} &= \begin{pmatrix}
        8.0 & 26.1\\ 26.1 & 145.3
    \end{pmatrix}
\end{align*}
we note an increased incentive for the MPC controller to bring the terminal state closer to rest when we insist to operate within the maximal operating range. To validate our results further, we let $x_0=\begin{pmatrix}
    7.99 & -1.27
\end{pmatrix}$,
which is just inside the lower right corner of the maximal control invariant set in Figure~\ref{fig:invsets}. Recall that the MPC controller that uses the invariant set of the LQR controller and the terminal cost defined by $P_{\infty}$ above requires a planning horizon of at least $26$ samples to operate from this initial value. In contrast, the simulations shown in Figure~\ref{fig:sim2d} show how the MPC controller drives the system to rest for all horizon lengths.

\begin{figure}[htb]
\centerline{\includegraphics[width=\hsize]{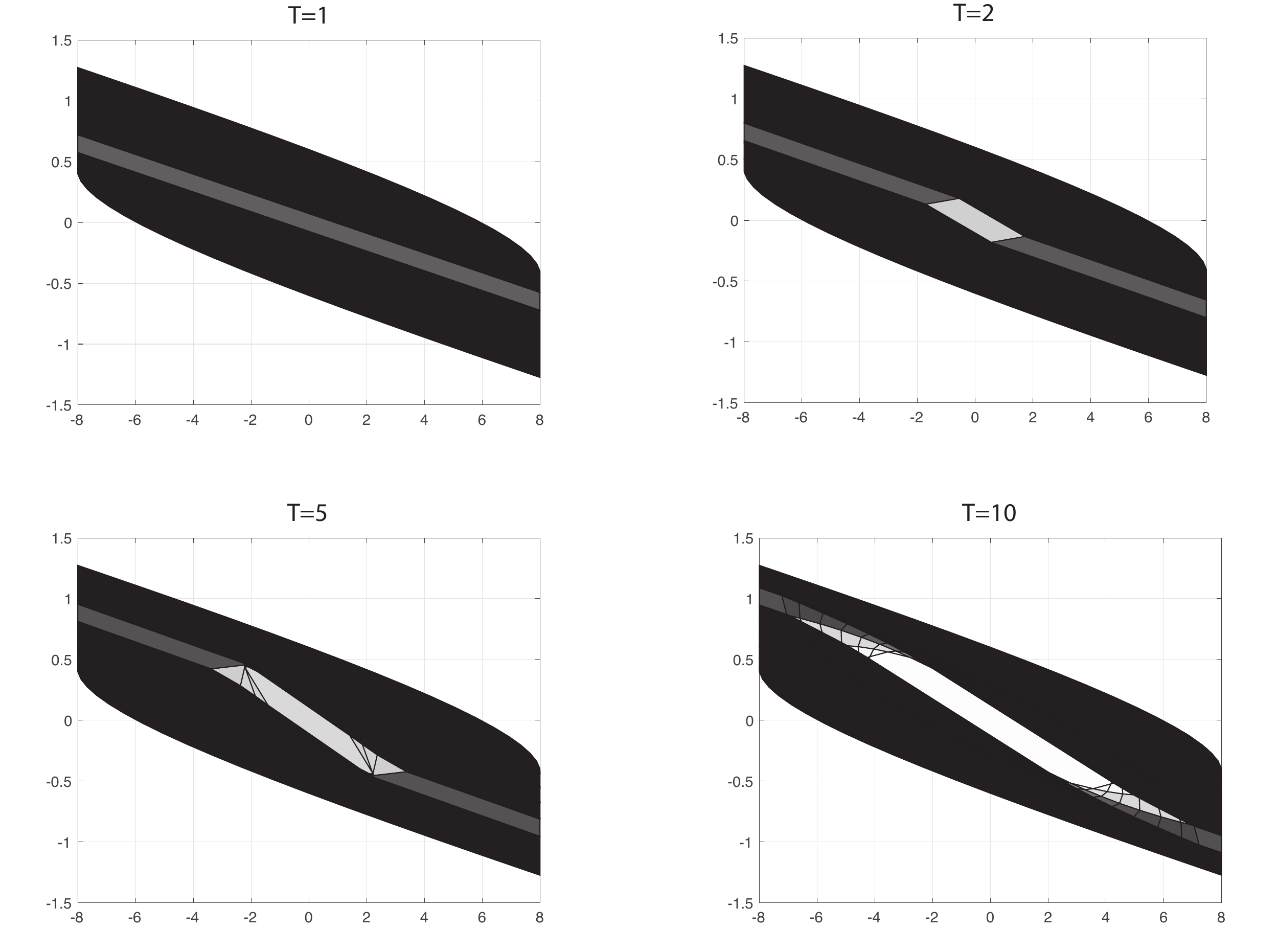}}
\caption{Partitions of explicit MPC policies for different $T$. Light colors indicate a local behavior close to $u_t=-L_{\infty}x_t$.}
\label{fig:xmpc_fbpart}
\end{figure}

Figure~\ref{fig:xmpc_fbpart} shows the partitions of the corresponding explicit MPC control laws. Here, a lighter color indicates that the local behavior in the region is closer the infinite-horizon LQR controller, with white for perfect agreement and black for  saturation. We notice how the linear behavior around $x=0$ becomes increasingly close to the infinite-horizon LQR controller as the prediction horizon increases.

\subsection{Stability in absence of terminal set contractivity}
Next, we consider the constrained linear system given by
\begin{align*}
    A &= \begin{pmatrix}
        0 & 1\\ -1 & 0
    \end{pmatrix},\qquad B=
    \begin{pmatrix}
        0\\ 1
    \end{pmatrix}
\end{align*}
and $\mathcal{X}=\{ x \;\vert\; \Vert x\Vert_{\infty}\leq 5\}$ and $\mathcal{U}=\{ u \;\vert\; \vert u\vert\leq 1 \}$. It turns out that the maximal control invariant set $\mathcal{C}_{\infty}$ is $\mathcal{X}$ itself. 

Since ${\mathcal C}_{\infty}$ is not contractive,
this system is challenging for approaches such as \cite{GrP:13} that use the Minkowski functional of the maximal control invariant set as terminal cost. To see that $\mathcal{C}_{\infty}$ is not contractive, note that  $x_t=(
    0\; 5
) \in \partial \mathcal{C}_{\infty}$. In addition, the first component of $x_{t+1}$ is equal to $5$ for all $u_t$, which means that also $x_{t+1}\in \partial \mathcal{C}_{\infty}$. Hence, the set is not contractive, and the corresponding Minkowski functional will not decrease for all $x\in {\mathcal C}_{\infty}$. 

However, the system admits a quadratic terminal cost that can be found by the approach developed in this paper. With $Q=I$ and $R=1$, the algorithm in Setion~\ref{sec:algorithm} returns 
\begin{align*}
    Q_T &= \begin{pmatrix}
    9.65  &  0.50\\
    0.50  & 10.67
    \end{pmatrix}
\end{align*}
which is the infinite-horizon cost-to-go for the feedback law $u_t=(0.1\, -0.1)x_t$. In fact, all state feedback laws on the form $u_t=-\begin{pmatrix}
    l_1 & l_2
\end{pmatrix}x$ with $\vert l_1+0.1\vert+\vert l_2\vert \leq 0.1$ and $l_1$ and $l_2$ not both equal to zero are admissible in ${\mathcal X}$, render the closed-loop system asymptotically stable and ${\mathcal X}$ invariant. 

\subsection{A higher-order system}
\begin{figure}[htb]
\centerline{\includegraphics[width=0.95\hsize]{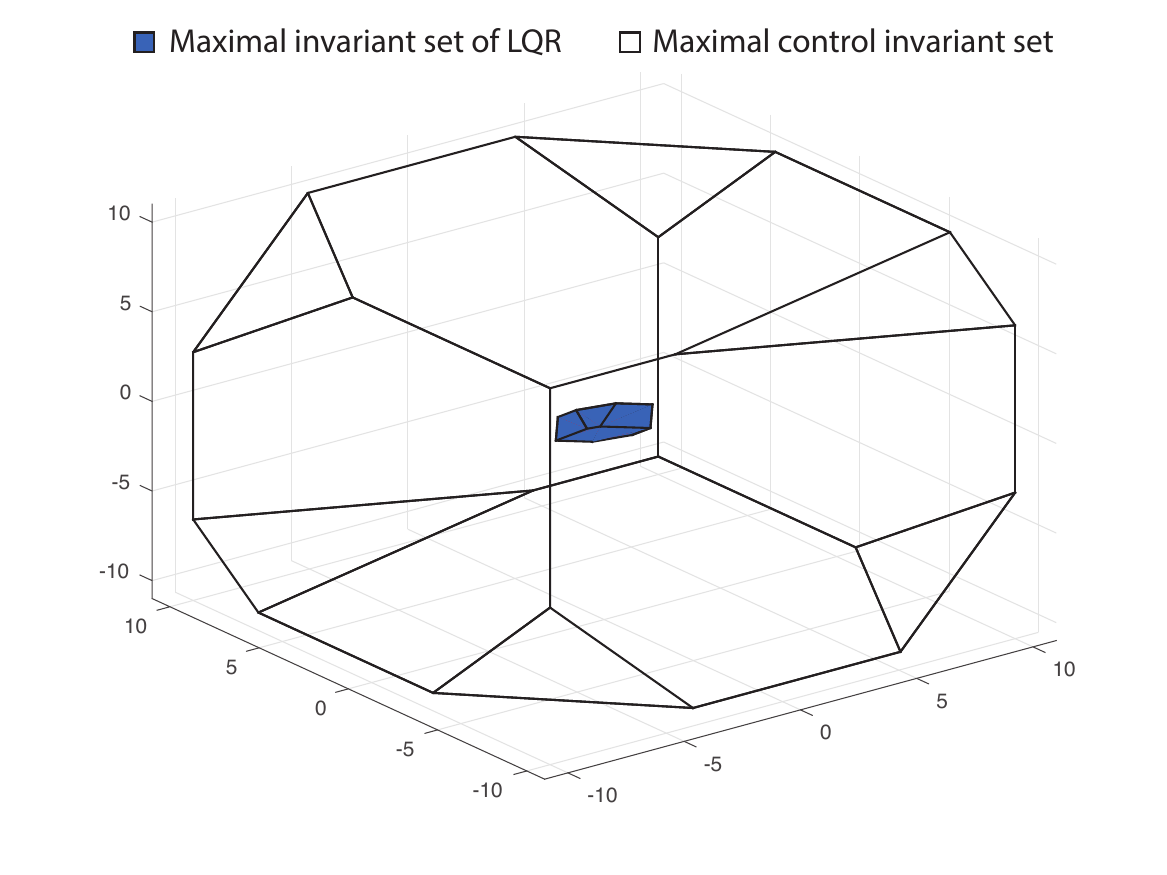}}
\caption{Maximal invariant set is much larger than the invariant set of the LQR controller.}
\label{fig:invsets_large}
\end{figure}

\begin{figure}[htb]
\centerline{\includegraphics[width=0.925\hsize]{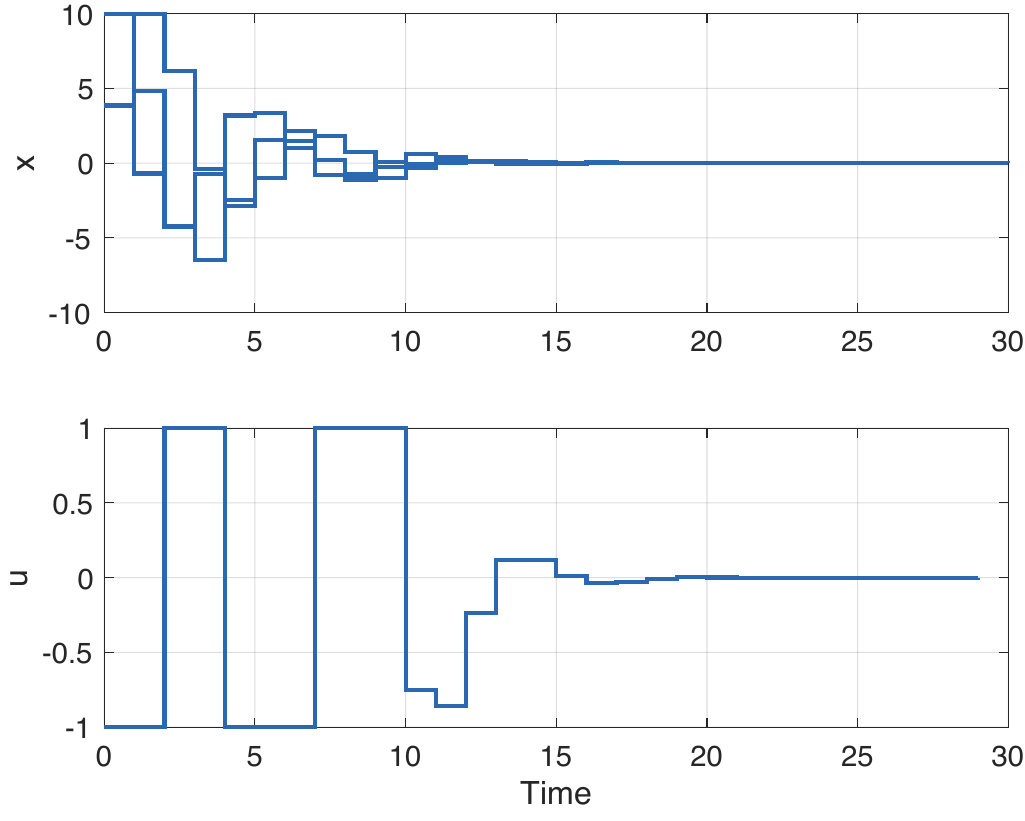}}
\caption{Initial value responses for MPC controller with terminal penalty and terminal cost defined by LQR solution and horizon $T=10$ (the minimal horizon required).}
\label{fig:large_sims}
\end{figure}

As another example, we consider the system
\begin{align*}
    x_{t+1} &= 
    \begin{pmatrix}
    0.48 &   0.45 &   0.38 \\
   -0.13 &   0.52 &  -0.54 \\
   -0.58 &   0.32 &   0.40    
    \end{pmatrix}x_t + 
    \begin{pmatrix}
        0.15 \\
         0.00 \\
    0.14 
    \end{pmatrix}u_t
\end{align*}
with cost given by $Q=10 I$ and $R=1$, and constraints $\mathcal{X}=\{ x \;\vert\; \Vert x\Vert_{\infty}\leq 10\}$ and $\mathcal{U}=\{ u \;\vert\; \Vert u\Vert_{\infty}\leq 1\}$. In this case, the maximal control invariant set is much larger (a factor 500x) in volume than the invariant set of the LQR controller, see Figure~\ref{fig:invsets_large}.
Yet, the quadratic upper bounds are quite similar: our procedure finds
\begin{align*}
    Q_T &=
    \begin{pmatrix}
    23.98 &  -1.69 &   1.03 \\
   -1.69  & 23.25 &   0.78 \\
    1.03  &  0.78 &  24.65
    \end{pmatrix}
\end{align*}
while the Riccati solution is 
\begin{align*}
    P_{\infty} &= 
    \begin{pmatrix}
   21.18  &  0.59 &    0.04\\
    0.59  & 18.21  & -1.95 \\
    0.04 &  -1.95 &   18.97
    \end{pmatrix}
\end{align*}
We simulate the closed-loop system from $x_0=\begin{pmatrix} -10 & 10 & 0\end{pmatrix}^{\top}$. The MPC controller based on the LQR invariant set and terminal cost requires a horizon of $T=10$ for this initial value, and yields the closed-loop response shown in Figure~\ref{fig:large_sims}. When we use the maximal control invariant set and our upper bound cost, on the other hand, we can use all horizon lengths and get practically indistinguishable responses even for $T=1$.

\section{Conclusions} \label{sec:conclusions}
We have presented a numerical procedure for computing a quadratic cost-to-go for linear MPC controllers that is valid for the maximal control invariant sets. This results in the largest possible set of recursively feasible states for the MPC controller, while the closed-loop stability follows from the standard stability proof for linear MPC.

We believe that the suggested procedure could be adapted to many scenarios beyond linear MPC and maximal control invariant sets. In essence, the approach only relies on our ability to compute \emph{a} control invariant set (not necessarily the maximal one), recover admissible vertex controls, and to interpolate these into a piecewise linear feedback law. We leave such extensions as future work. 

\bibliographystyle{unsrt}
\bibliography{maxinv}
\end{document}